\documentclass[reqno,11pt]{amsart}
\usepackage{amsmath, latexsym, amsfonts, amssymb, amsthm, amscd}
\usepackage{graphics,epsf,psfrag}
\numberwithin{equation}{section}

\newtheorem{theorem}{Theorem}[section]

\newtheorem{proposition}[theorem]{Proposition}

\newtheorem{rem}[theorem]{Remark}

\DeclareMathOperator{\sign}{\mathrm{sign}}

\newcommand{\ind}{\mathbf{1}}

\newcommand{\R}{\mathbb{R}}
\newcommand{\Z}{\mathbb{Z}}
\newcommand{\N}{\mathbb{N}}

\renewcommand{\hat}{\widehat}

\newcommand{\cN}{{\ensuremath{\mathcal N}} }

\newcommand{\bP}{{\ensuremath{\mathbf P}} }
\newcommand{\bE}{{\ensuremath{\mathbf E}} }

\DeclareMathSymbol{\leqslant}{\mathalpha}{AMSa}{"36} 
\DeclareMathSymbol{\geqslant}{\mathalpha}{AMSa}{"3E} 
\DeclareMathSymbol{\eset}{\mathalpha}{AMSb}{"3F}     
\newcommand{\dd}{\,\text{\rm d}}             

\newcommand{\bbC}{{\ensuremath{\mathbb C}} }

\newcommand{\bbP}{{\ensuremath{\mathbb P}} }

\newcommand{\ga}{\alpha}
\newcommand{\gb}{\beta}
\newcommand{\gd}{\delta}
\newcommand{\gep}{\varepsilon}       

\newcommand{\gG}{\Gamma}

\newcommand{\gD}{\Delta}

\makeatletter
\def\captionfont@{\footnotesize}
\def\captionheadfont@{\scshape}

\long\def\@makecaption#1#2{%
  \vspace{2mm}
  \setbox\@tempboxa\vbox{\color@setgroup
    \advance\hsize-6pc\noindent
    \captionfont@\captionheadfont@#1\@xp\@ifnotempty\@xp
        {\@cdr#2\@nil}{.\captionfont@\upshape\enspace#2}%
    \unskip\kern-6pc\par
    \global\setbox\@ne\lastbox\color@endgroup}%
  \ifhbox\@ne 
    \setbox\@ne\hbox{\unhbox\@ne\unskip\unskip\unpenalty\unkern}%
  \fi
  \ifdim\wd\@tempboxa=\z@ 
    \setbox\@ne\hbox to\columnwidth{\hss\kern-6pc\box\@ne\hss}%
  \else 
    \setbox\@ne\vbox{\unvbox\@tempboxa\parskip\z@skip
        \noindent\unhbox\@ne\advance\hsize-6pc\par}%
\fi
  \ifnum\@tempcnta<64 
    \addvspace\abovecaptionskip
    \moveright 3pc\box\@ne
  \else 
    \moveright 3pc\box\@ne
    \nobreak
    \vskip\belowcaptionskip
  \fi
\relax
}
\makeatother
\def\writefig#1 #2 #3 {\rlap{\kern #1 truecm
\raise #2 truecm \hbox{#3}}}

\newcommand{\Kbar}{{\overline{K}}}

\begin{document}

\title[Critical pinning and renewal convergence rates]{
Renewal convergence rates
and correlation decay for homogeneous pinning models 
}

\author{Giambattista Giacomin}
\address{Universit{\'e} Paris 7 -- Denis Diderot
and 
 Laboratoire de Probabilit{\'e}s et Mod\`eles Al\'eatoires (CNRS U.M.R. 7599), 
 U.F.R. Mathematiques, Case 7012, 2 place Jussieu, 75251 Paris cedex 05, France}
\email{giacomin\@@math.jussieu.fr}

\date{\today}

\begin{abstract}
A class of discrete renewal processes with super-exponentially decaying inter-arrival distributions coincides with the  infinite volume limit of general homogeneous
  {\sl pinning} models in their localized phase. 
  Pinning models are statistical mechanics systems 
to which a lot of attention has been devoted both for their relevance
for applications and    
 because  they are {\sl solvable} models exhibiting 
  a non-trivial phase transition. 
  The spatial decay of correlations
  in these systems is   directly  mapped to 
  the speed of convergence to equilibrium for the associated renewal processes.
 We show that close to criticality, under general assumptions, 
   the correlation decay rate, or the renewal convergence rate,
     coincides with the inter-arrival decay rate. We also show  that, in general, 
     this is false away from criticality.
 Under a stronger assumption on the inter-arrival distribution we
 establish a local limit theorem, capturing thus the sharp
 asymptotic behavior of correlations.
\end{abstract}

\maketitle

\section{Introduction and main results}
\label{sec:intro}

\subsection{Set--up and generalities}
\label{sec:Kb}

We consider the probability discrete density $K(\cdot)$
concentrated on
 $\N=\{1,2, \ldots\}$. We choose $K(\cdot)$ such that
 for some $\ga >0$ and some  function $L(\cdot)$ which is 
 slowly varying at infinity we have
 \begin{equation}
 \label{eq:weakHyp}
 \Kbar (N) \, :=\, 
 \sum_{n> N}K(n)\stackrel{N \to \infty} \sim \frac{L(N)}{\ga N^\ga},
 \end{equation}
 where we are using the notation $a_N \stackrel{N \to \infty} \sim b_N$
 when $\lim_{N \to \infty} a_N/b_N=1$. We assume that
 $K(\cdot)$ is aperiodic, that is 
 $\text{gcd}\{n: \,  K(n)>0\}=1$. 
 We recall that 
 a function $L(\cdot)$ defined on the positive semi-axis is slowly varying at infinity
 if it is positive, measurable  and if 
  $\lim_{t \to \infty}L(ct)/L(t) =1$ for every
 $c>0$. We refer to \cite{cf:RegVar} for the full theory of slowly varying functions,
 recalling simply that both $L(t)$ and $1/L(t)$ are much smaller than $t^\gd$ (as $t \to \infty$),
 and this for 
 any $\gd >0$.
 
We point out that \eqref{eq:weakHyp} and aperiodicity are implied by
\begin{equation}
 \label{eq:strongHyp}
 K(n)\stackrel{n \to \infty} \sim \frac{L(n)}{n^{1+\ga}}.
 \end{equation}
\medskip

Starting from $K(\cdot)$, we introduce a family of  discrete probability densities
indexed by
 $b\ge 0$:
\begin{equation}
\label{eq:Kb}
K_b(n)\, :=\, c(b) K(n) \exp(-bn),
\end{equation}
and $c(b)= 1/ \sum_n K(n) \exp(-bn)$ (of course $c(0)=1$).

Our attention focuses on   
the  renewal process $\tau (b):=\{\tau_0(b), \tau_1(b), \tau_2(b), \ldots\}$ 
with inter-arrival law $K_b(\cdot)$, that is  the process defined by $\tau_0(b)=0$ and by the requirement that
$\{ \tau_{j+1}(b)-\tau_j(b)\}_{j=0, 1, \ldots}$ is a sequence of IID random variables 
 and $\bP(\tau_1(b)=n)=K_b(n)$. Note that $\tau (b)$ is an increasing sequence
 of almost surely finite numbers and it 
 can be looked upon equivalently as a sequence of random variables (in fact,
 a random walk with positive increments) or as a random subset of $\N \cup \{ 0\}$.
 With this second interpretation we introduce the so called
 {\sl mass renewal} function, that is
 \begin{equation}
 u_b(n) \, := \, \bP \left( n \in \tau (b)\right),
 \end{equation} 
 so that $u_b(n)$ is the probability that the site $n$ is visited 
 by the renewal. Note that  $u_b(0)=1$ and, since $K(\cdot)$
 is   aperiodic, there
 exists $n_0>0$ such that  $u_b(n)>0$
 for every $n\ge n_0$.
  
 \subsection{The Renewal Theorem and refinements to it}
We now make an excursus in the general renewal theory on the integer numbers.
We consider thus a general
renewal process with $\tau_0=0$ and with inter-arrival taking values in $\N$. 
For this we introduce the notation
$F(n):= \bP (\tau_1=n)$, while the mass renewal function is  denoted 
by $u(\cdot)$.
The classical Renewal Theorem (see {\sl e.g. } \cite{cf:asmussen})
says that, if $F(\cdot)$ is aperiodic, we have
 \begin{equation}
 \label{eq:RenTh}
u(\infty) \, :=\, \lim_{n \to \infty} u(n) \, =\, \frac 1{\bE[\tau_1]}\in [0, 1].
 \end{equation}
\smallskip

Much effort has been put into refining such a result. 
Refinements are of course a very natural question when
 $\bE[\tau_1]=+\infty$ ({\sl e.g.} \cite{cf:Doney,cf:GL}), 
 as well as  if $\bE[\tau_1]<+\infty$. In the latter case sharp estimates on
 $u (n) -u(\infty)$ have been obtained for sub-exponential tail decay
 of the inter-arrival distribution, like for example in the case of
 $F(\cdot)=K(\cdot)$ and $K(\cdot)$ as
 in \eqref{eq:strongHyp} (we refer to \cite{cf:Grubel} and references therein). 

\smallskip 
 When instead
 the inter-arrival distribution decays super-exponentially, like for example if 
 $F(\cdot)=K_b(\cdot)$ with $b>0$, 
 general sharp results are harder to obtain. What can be proven in general
 in fact is that,
 if there exists $c_1>0$ such that $\lim_{n \to \infty }\exp(c_1 n) F (n)=0$, 
then there exists $c_2>0$ such that $\lim_{n \to \infty }\exp(c_2 n)\vert u(n) -u(\infty) \vert=0$.
However the precise decay, or even only the exponential asymptotic behavior (that is 
the supremum of the  values of $c_2$ for which the previous equality holds), 
in general does not depend only on the tail behavior of the inter-arrival probability.
This is definitely  a very classical problem
\cite{cf:Kendall,cf:Heathcote}, and a number of results have been proven in
specific instances (see {\sl e.g} \cite{cf:BL,cf:MT}). 
We are now going to 
treat this point in some detail.

\subsection{ On super-exponentially decaying
inter-arrival laws}

 From the very definition of renewal process one
directly derives
the equivalent expressions
\begin{equation}
\label{eq:rec}
u(n)\, =\, \ind _{\{0\}} (n) + \sum_{j=0}^{n-1} u(j)F (n-j)
\ \ \ \ \  \text{ and }\ \ \  \ \   \hat u(z) \, =\, \frac 1{1- \hat F (z)},
\end{equation}
with the notation $\hat f (z) = \sum_{n=0}^\infty z^n f(n)$ ($\hat f(\cdot)$ is the $z$-{\sl transform}
of $f(\cdot)$) and 
$z$ is a complex number.
Of course $\hat f(\cdot)$ is a power series and $\vert z\vert$  a priori has to be chosen
smaller than the radius of convergence,
 which, for the two series appearing in \eqref{eq:rec}, is at least  $1$.

As a matter of fact, 
we are interested (in particular) in the radius of convergence
of the series
\begin{equation}
\label{eq:BL}
\gD (z)\, :=\, \sum_{n=0}^\infty (u(n)-u(\infty))z^n
 \, =\, \frac1{1-\hat F (z)} - \frac 1{\bE[\tau_1] (1-z)}
.
\end{equation}
If we assume that $\limsup_{n\to \infty} \exp (cn) F (n) < \infty$ for some $c>0$,
the radius of convergence of $\hat F( \cdot)$ is at least $\exp(c)$, however it is not at all clear 
 that the radius of convergence of $\gD(\cdot)$ coincides 
 with the radius of convergence of $\hat F (\cdot)$. In reality the problem
 does not come from the singularity at $z=1$ ($\hat F (1)=1$) since it is 
 easily seen that it is removable. Notice also that, when $F(\cdot)$ is aperiodic,
 $\hat F(z)=1$ on the unit circle only if $z=1$. 
 However there may be other
 solutions $z$ to $\hat F ( z)=1$ for $z$ within the
 radius of convergence of $\hat F (\cdot)$.   
 And it may even happen that  
 $\gD (\cdot)$ can be analytically continued 
 beyond the radius of convergence of $\hat F (\cdot)$.  
 Let us make this clear by giving two explicit examples:
 \begin{itemize}
 \item $F (1)=1-p$, $F(2)=p$ and
 $    F (n)=0$ for $n=3,4, \ldots$ ($p\in (0,1)$). The radius of convergence of 
 $\hat F (\cdot)$ is $\infty$, but $\gD (z)=p/((1+p)(1+pz))$ and therefore
  the radius of convergence of $\gD(\cdot)$ is $1/p$, and in fact, by expanding 
  $\gD (z)$ around $z=0$, we obtain 
  $u(n)-u(\infty)= (-p)^n (p/(1+p))$ for  $n=1,2, \ldots$. 
 \item $F(n) = p^n (1-p)/p$, $p \in (0,1)$. In this case
 the radius of convergences of  $\hat F (\cdot)$ is $1/p$,
 but $\gD(z)=p$ for every $z$, so the radius of convergence
 is $\infty$ and in fact $u (n)-u(\infty)=0$ for every $n\ge 1$. 
 \end{itemize}
 
 These examples show  that 
 the tail decay of $u(\cdot)-u(\infty)$   
  may have little to do with the tail decay of
   the $F(\cdot)$: in particular, changing fine details of $F(\cdot)$ may have a drastic
   effect on the decay  of $u(\cdot)-u(\infty)$.
 For further examples of such a  behavior 
  see in particular \cite{cf:BL}, but also Section~\ref{sec:examples} below.
 
 The main purpose of this note is, however, to point out that, in 
 a suitable class of renewal processes, motivated by 
 statistical mechanics modeling (see Subsection~\ref{sec:pinning}), 
 the tail decay of $u(\cdot)-u(\infty)$   
 is closely linked  with the tail decay of
   the $F(\cdot)$. We are in fact  going to show
 that if $F(\cdot)=K_b (\cdot)$, that is in the set-up of
 \S~\ref{sec:Kb}, the decay rate of $\{u_b(n)-u_b(\infty)\}_n$
  is equal  to the decay rate  of $K_b(\cdot)$,
 if $b$ is sufficiently small. And under the stronger
 hypothesis \eqref{eq:strongHyp} we control the sharp asymptotic 
 behavior of $u_b(n)-u_b(\infty)$.

 \subsection{Main result}
 With the set-up of \S~\ref{sec:Kb} we have the following:
 \medskip
 
 \begin{theorem}
 \label{th:main}
  Given $K(\cdot)$ call $b_0(\in[0, \infty])$
 the infimum of the values of $b>0$ such that there exists $z$
 satisfying $1< \vert z\vert \le \exp(b)$ and 
  $ \hat{K_b} (z)=1$.
 \begin{enumerate} 
 \item
 For every choice of $K(\cdot)$ satisfying
 \eqref{eq:weakHyp} we have 
 $b_0 \in (0, \infty]$ and  for every  $b \in (0, b_0]$ we have 
 \begin{equation}
 \label{eq:main0}
 \limsup_{n\to \infty} \frac 1n \log \left \vert u_b(n)-u_b(\infty)
 \right\vert \, =\, -b,
 \end{equation} 
 while for $b>b_0$ we have 
 \begin{equation}
 \label{eq:main0.1}
 \limsup_{n\to \infty} \frac 1n \log \left \vert u_b(n)-u_b(\infty)
 \right\vert \, \ge \, -b.
 \end{equation} 
 \item
 For every choice of $K(\cdot)$ satisfying
 \eqref{eq:strongHyp}
 we have 
 that   for every  $b \in (0, b_0)$ 
 \begin{equation}
 \label{eq:main1}
 u_b(n)-u_b(\infty)\stackrel{n \to \infty} \sim \frac{K_b(n)}{(c(b)-1)^2}
  ,
 \end{equation}
 which implies 
 \begin{equation}
 \label{eq:main1.1}
 \lim_{n\to \infty} \frac 1n \log \left( u_b(n)-u_b(\infty)
 \right) \, =\, -b.
 \end{equation} 
 \end{enumerate}
 \end{theorem}  
 \medskip
 
 \begin{rem} \rm
 When there exists $z_0$, $1< \vert z_0\vert < \exp(b)$, such that
 $\hat{K_b}(z_0) =1$ 
 (therefore $b>b_0$) one can easily write down the sharp asymptotic
 behavior of $\{ u_b(n)-u_b(\infty)\}_n$ in terms of the values of $z_0$
 with minimal $\vert z_0\vert$. As a matter of fact one has
 \begin{equation}
 \label{eq:main0.1.2}
 \limsup_{n\to \infty} \frac 1n \log \left \vert u_b(n)-u_b(\infty)
 \right\vert \, = \, - \log \vert z_0\vert \, >\, -b ,
 \end{equation} 
 but
 the sequence   changes sign
 infinitely often and, in general, the superior limit cannot be replaced by a
 limit (see Section~\ref{sec:examples} for details).
 In  Section~\ref{sec:examples} we also provide explicit examples showing that
 $b_0$ can be arbitrarily small by 
 choosing $K(\cdot)$ suitably. 
 In all the examples we have worked out 
 the inequality in \eqref{eq:main0.1} is strict (for every $b>b_0$),
 but it is unclear to us whether or not this is a general phenomenon.
\end{rem}
 \medskip

The proof of Theorem~\ref{th:main}(1) can be found in 
Section~\ref{sec:radius}
which
is devoted to the study of $R_b:=1/ \limsup_n \vert u_b(n)-u_b(\infty)\vert ^{1/n}$,
which of course is the radius of convergence of
$\gD_b (\cdot)$, and to establishing that $b_0$ is not zero.
Theorem~\ref{th:main}(2) follows instead by a direct application
of a well established technique
 \cite{cf:CNW}: we detail this application in Section~\ref{sec:sharp}.
 We point out that
the validity of the results in \cite{cf:CNW} go beyond
the assumption 
\eqref{eq:strongHyp}, but we do make use of the regularly varying
character of $K(\cdot)$ in establishing $b_0>0$.
A closer look at the proof of $b_0>0$
however 
shows that when $\sum_n nK(n)< \infty$ ({\sl cf.} \eqref{eq:Rieman2.1.2})
the regular variation
property is used only marginally and in fact Theorem~\ref{th:main}
holds also for a number of sub-exponential ({\sl c.f.} \cite{cf:RegVar})
distributions $K(\cdot)$ beyond our assumptions. For example
Theorem~\ref{th:main}
holds also for  $K(n)=L(n)n^q\exp(-n^\gamma)$,
with $q\in \R$ and $\gamma \in (0,1)$. 

 .

 \subsection{Homogeneous pinning models and decay of correlations}
 \label{sec:pinning}

What motivated,  and what even suggested the validity of the results in this note,
is the behavior near criticality of homogeneous pinning models.
As it as been pointed out in particular in \cite{cf:Fisher}, a 
 large class of physical models  boils down to 
 a class of Gibbs measures that, in mathematical terms, are just 
 obtained from discrete renewal processes {\sl modified} by introducing
 an exponential   weight, or Boltzmann factor, depending 
 on $\cN_N(\tau):=\vert \tau \cap (0, N]\vert$. More precisely
 if $\bbP$ is the law of $\tau$ and the latter is the renewal 
 sequence with inter-arrival  distribution $K(\cdot)$, we consider
 the family of probability measures  $\left\{ \bbP_{N, \gb}\right\}_{N \in \N}$ defined
 by 
 \begin{equation}
 \label{eq:RN}
 \frac{\dd \bbP_{N, \gb}}{\dd \bbP} (\tau)\, =\,
 \frac 1{Z_{N ,\gb}}
 \exp\left(
 \gb \cN_N(\tau)
 \right),
 \end{equation}
 with $Z_{N , \gb}$ the normalization constant.
 Then one can show (\cite{cf:CGZ},\cite[Ch.~2]{cf:G})
  that the weak limit $\bbP_{\infty,\gb}$ of 
 $\left\{ \bbP_{N, \gb}\right\}_{N \in \N}$
 exists for every $\gb\in \R$ (to be precise, this statement holds for every $\gb$
 assuming \eqref{eq:strongHyp}, but it holds also assuming only
 \eqref{eq:weakHyp} if $\gb>0$). The parameter $\gb$ actually plays
 a crucial role. In fact if $\gb< 0$  then $\tau$, under
 $\bbP_{\infty,\gb}$, is a transient renewal and it contains 
therefore only a finite number of points  (this is the so-called {\sl delocalized phase}).
If instead $\gb>0$ then $\tau$, again  under
 $\bbP_{\infty,\gb}$, is a positive recurrent renewal
 with inter-arrival distribution given by $K_b(\cdot)$,
 with $b=b(\gb)$  unique real solution of $\sum_n K(n) \exp(-bn)=\exp(-\gb)$
 (this is the {\sl localized phase}).
 Note that if $\gb\searrow 0$, then $b\searrow 0$.
 
 We point also out that it is not difficult  to see that $b$ 
 coincides with the limit as $N $ tends to infinity of
 $(\log Z_{N , \gb})/N$ and it is hence the {\sl free energy}
 of the system \cite[Ch.~1]{cf:G}.
 In  \cite{cf:Fisher} and, more completely in \cite[Ch.~2]{cf:G},
 one can find 
the analysis of $b(\gb)$ as $\gb \searrow 0$.
 
 \smallskip
 
 As a consequence $\tau(b)$, for $b>0$, does describe 
 the localized regime of an infinite volume statistical mechanics system: 
 if $b$ is small, the system is close to {\sl criticality}.
 The {\sl correlation length} is a 
  key quantity in statistical mechanics, see {\sl e.g.}  \cite{cf:Fisher}. Moreover
 it is expected to scale {\sl nicely}
  with $\gb$ (or, which is equivalent,
 with $b$) approaching criticality, 
 typically as $\gb$ to some (negative) power, possibly times
 {\sl logarithmic} corrections.
 The correlation length may be defined 
 by introducing first the correlation function:
 \begin{equation}
 \begin{split}
 \mathbf{c} (n)\, &:= \, \lim_{m \to \infty}
 \frac{
 \bP\left( m\in \tau (b), \, m+n \in \tau (b) \right) 
 -  \bP\left( m\in \tau (b)\right) \bP\left(  m+n \in \tau (b) \right)
 }{
 \sqrt{ \bP\left( m\in \tau (b)\right)
\left(1- \bP\left( m\in \tau (b)\right) \right)
 \bP\left( 
  m+n \in \tau (b) \right)
  \left(1- \bP\left( m+n\in \tau (b)\right) \right)
  } } 
  \\
  &  \phantom{:}=\,
  \frac{\bE\left[ \tau_1(b) \right]}{\bE\left[ \tau_1(b) \right]-1} 
  \left(
 \bP\left( n\in \tau (b) \right)  -\frac 1{\bE\left[ \tau_1(b) \right]}\right),
 \end{split}
 \end{equation}
 where we have used the Renewal Theorem.
Then  the correlation length is just one over the decay rate $\xi(b)$ of $\mathbf{c} (\cdot)$:
$\xi (b) := -1/\limsup_{n \to \infty } n^{-1}\log \vert \mathbf{c} (n)\vert $ and therefore 
\begin{equation}
\xi (b)\, =\, -1/\limsup_{n \to \infty } n^{-1}\log \vert u_b (n)-u_b(\infty)\vert,
\end{equation} 
so that Theorem~\ref{th:main} guaranties that 
\begin{equation}
\label{eq:xi}
 \xi (b)  \stackrel{b \searrow 0}\sim \frac 1 b,
\end{equation}
which roughly can be rephrased 
by saying that
the correlation length, close to criticality, scales
like one over the free energy.

\medskip

On physical grounds \eqref{eq:xi}, or rather the weaker form $
\log \xi(b) \sim -\log b$,
  is certainly expected \cite{cf:Fisher}.
A proof of \eqref{eq:xi} has been given 
in \cite{cf:Fabio1} by coupling arguments for the case in which
$K(\cdot)$ is given by the return times of a simple random walk
(and the proof is given also for disordered models).
The result actually holds as an equality for every $b$ (like the
case presented in \S~\ref{sec:basic} below: we point out 
that for $\ga=1/2$ the distribution
$K(\cdot)$ treated in  \S~\ref{sec:basic}
coincides with the distribution of the returns to zero
of a simple random walk in the sense that $K(n) $ is the probability
that the first return to zero of a simple random walk happens 
at time $2n$).
In general coupling arguments yield sharp results on the rate
when suitable monotonicity properties 
are present (see in particular \cite{cf:LT}): the returns of a simple random walk
are in this class.
In absence of monotonicity properties coupling
arguments usually yield only upper bounds on the speed
of convergence (and hence lower bounds on the rate, see 
\cite{cf:asmussen} and references therein):
in \cite{cf:Fabio2} a coupling argument is given for disordered pinning models 
and it yields in our homogeneous set-up
that $\limsup_{b \searrow 0} \log \xi(b)/\log(b) \le -1$,
under the stronger hypothesis \eqref{eq:strongHyp}.

We conclude this introduction with two important remarks:

\begin{rem}\rm
Some 
of the papers we have referred to  (in particular  \cite{cf:BL,cf:MT})
aim at 
 explicit bounds  that hold for every $n$, possibly at the
expense of sharp asymptotic results. 
Also in our set-up
the question of obtaining more quantitative estimates, particularly
when $b \searrow 0$, is important and relevant for the applications.
\end{rem}

\begin{rem}\rm
The class of pinning models we have considered contains  the so called
$(1+d)$--{\sl dimensional} pinning models. 
The name comes from the {\sl directed} viewpoint
on Markov chains: if one considers a Markov chain $S$ with state    space
$\Z ^d$, the state space of the {\sl directed process} $\{ (n, S_n)\}_n$ is 
$\Z ^{1+d}$. The renewal structure in this case is simply given
by the successive returns to $0\in \Z^d$ by $S$ or,
equivalently, by the intersections of the directed process with
the line $\{(n,0) \in \Z^{1+d}: \, n=0,1, 2 , \ldots\}$.  
This viewpoint is important in order to understand the spectrum
of applications of pinning models. We are not going to discuss 
this further here, and we refer to \cite{cf:G,cf:Yvan}, but we do 
point out that precise estimates catching the order
of magnitude of  the correlation length
in a class of $(d+1)$--dimensional pinning models, 
{\sl i.e.} Gaussian effective 
surfaces in a $(d+1)$--dimensional space pinned at 
an hyper-plane, have been obtained in     \cite{cf:BV}.
\end{rem}

\section{The radius of convergence of $\gD_b(\cdot)$}
\label{sec:radius}

In this section we work
in the most general set-up,
{\sl i.e.} 
we assume \eqref{eq:weakHyp}.
Recall the definition of $b_0$ from the statement of Theorem~\ref{th:main}.

\medskip

\begin{proposition}
\label{th:Rb}
$R_b \le \exp(b)$ and, for every choice of $K(\cdot)$, $b_0>0$ and therefore
 $R_b= \exp(b)$
for $b \in (0, b_0]$.
\end{proposition} 
\medskip

Note that this result
implies \eqref{eq:main0} and
 \eqref{eq:main0.1}.
\smallskip

\noindent
{\it Proof.}
We are going to show that $R_b\le \exp(b)$
by  making use only of $\hat {K_b}(\exp(b))<\infty$
and of the fact that the radius of convergence of
$\hat{K_b}(\cdot)$ is $\exp(b)$.

\smallskip
Of course we may assume that $\Delta_b(\cdot)$ is analytic
in the centered ball of radius $\exp(b)$, since otherwise there is nothing to prove.
Let us suppose that $\gD_b(\cdot)$ has an analytic extension to
the open ball of radius $R>\exp(b)$.
From \eqref{eq:BL} we immediately derive an expression for
$\hat{K_b}(z)$ in terms of $\Delta_b(z)$, for $\vert z \vert < \exp(b)$,
and this gives the  meromorphic extension 
of $\hat{K_b}(\cdot)$ to the centered ball of radius $R$.
However we know that the radius of convergence of $\hat{K_b}(\cdot)$
is $\exp(b)$ and that $\vert \hat{K_b}(z) \vert \le \sum _ n K (n) < \infty$
if $\vert z \vert =\exp(b)$. So the singularity of $\hat{K_b} (\cdot)$ cannot
be a pole and therefore  $\hat{K_b}(\cdot)$ does not have a meromorphic extension. 
This implies that $\gD_b(\cdot)$ cannot be  analytically continued beyond
the centered ball of radius $\exp(b)$.
\medskip

The question that we have to address in order to complete the proof
of Proposition~\ref{th:Rb}, that is proving $b_0>0$,  can be rephrased as:
 do there exist two sequences
$\{b_j\}_j$, $b_j \searrow 0$ and $\{z_j\}_j$, $1<\vert z_j\vert \le \exp(b_j)$
such that $\hat{K_b} (z_j)=1$ for every $j$? 
Of course, if this is not the case, $\hat{K_b }(z)\neq 1$
if $\log \vert z\vert(>0)$ is sufficiently small. 

We make some preliminary observations: first,
we may assume $\Im (z_j )\ge 0$, since
if $\hat{K_b}(z)=1$, we have $\hat{K_b}(\overline{z})=1$ too.
Then let us remark that, by writing $z_j=r_j \exp(i\theta_j)$, 
we can pass to the limit in the equation $\hat{K_{b_j}} (z_j)=1$:
by the Lebesgue Dominated Convergence 
Theorem we have that every limit point $(1,\theta)$ of $\{(r_j , \theta_j)\}_j$
satisfies
\begin{equation}
\sum_n K(n) \exp(i n \theta)=1,
\end{equation}
which gives $\theta=0$ by aperiodicity.
This tells us that, for $b$ small, singularities have necessarily
positive real part and small imaginary part (in short, they are close 
to $1$). Moreover, by monotonicity, 
we see that the imaginary part cannot be zero
(and therefore we  assume that it is positive, since solutions
come in conjugate pairs).
\smallskip

Let us now assume by contradiction that there exists
a triplet of sequences 
\begin{equation}
\big(\{b_j\}_j,\, \{\gd_j\}_j, \, \{\theta_j\}_j \big),
\end{equation}
tending to zero, with the requirements that
$0\ge  \gd_j< b_j$, $\theta_j >0$ for every $j$ and
such that $\hat{K_b} (\exp(b_j-\gd_j)\exp(i\theta_j) ) =1$
for every $j$. Of course the triplet corresponds
to the poles of the associated $\gD_{b_j}(\cdot)$ function
at $z_j=\exp\left((b_j-\gd_j) +i \theta_j\right)$.
We are going to show that such a triplet does not exist
since we are able to extract
subsequences such that
\begin{equation}
\label{eq:aim}
\hat{K_b} (\exp(b_j-\gd_j)\exp(i\theta_j) ) \neq 1,
\end{equation}
for  every $j$ in the subsequence.

Let us consider the auxiliary sequence of non-negative numbers 
$\{ \gd_j / \theta_j \}_j$. By choosing a subsequence
we may assume that this sequence converges to a limit
point $\gamma \in [0, \infty]$. 

\smallskip

We consider first the case of $\ga\in (0,1)$. We distinguish the 
two cases $\gamma < \infty$ and $\gamma=\infty$.

\medskip

If $\gamma < \infty$ we have the asymptotic relation
\begin{equation}
\label{eq:Rieman2.1}
\sum_n K(n) 
\exp(-\gd_j n) \sin (\theta_j n)
\stackrel{j \to \infty}\sim
\theta_j ^\ga L(1/\theta_j) \int_0^\infty \frac{\exp(-\gamma s)\sin (s)}{s^{1+\ga}} \dd s\, ,
\end{equation}
that follows from
a Riemann sum approximation and the uniform convergence
property of slowly varying functions \cite[\S~1.5]{cf:RegVar} if the sum is restricted to
$\theta_j n\in (\gep , 1/\gep)$. The rest is then controlled for small $n$'s ($n \le \gep/\theta_j$) 
by replacing $ \sin (x)$ with $x$ and using  summation by parts which tells us that
$\sum_{n=1}^N n K(n)$ is equal to
$ \sum_{n=0}^{N-1} \Kbar (n) - N \Kbar (N)$ and the latter
 behaves for large values 
of $N$
as $N^{1- \ga}L(N)/(1-\ga)$ \cite[\S~1.5]{cf:RegVar}. For large $n$'s the rest is
controlled by using $\vert \exp(-\gd_j n) \sin (\theta_j n)\vert\le 1$.
Overall the absolute value of the rest is bounded by $c \theta_j^\ga L(1/\theta_j)
(\gep^{1-\ga}+ \gep^\ga)$ for some $c>0$, with $c$ not depending on $\gep$,
for $j $ sufficiently large (for example, $\theta_j < \gep$) and 
\eqref{eq:Rieman2.1} follows.

Observe that the left-hand side of \eqref{eq:Rieman2.1}
is asymptotically equivalent to
the imaginary part of $\hat{K_b} (\exp(b_j-\gd_j)\exp(i\theta_j) )$, apart for the
multiplicative  constant $c(b_j)= 1+o(1)\in \R$. 
The integral can be explicitly computed and it is equal to
\begin{equation}
\left( 1+\gamma^2\right)^{\ga/2} \gG(1-\ga) \sin\left( \ga \arctan (1/\gamma)\right),  
\end{equation}
which is positive for every $\gamma \in [0,\infty)$,
therefore for $j$ sufficiently large \eqref{eq:aim} holds 
(the definition of $\Gamma (\cdot)$ is recalled in Section~\ref{sec:examples}).

\medskip

If $\gamma = \infty $ instead we write
\begin{equation}
\label{eq:splitnew}
\sum_n K(n)
\exp(-\gd_j n) \sin (\theta_j n) \, =\, R_j^< + R_j^>,
\end{equation}
with $R_j^<$ the sum for $n \le \gep/\theta_j$ and $R_j^>$ is the rest
($0<\gep\le \pi/2$ is a fixed positive constant).
Setting $s_\gep := \sin(\gep)/\gep$ we have
\begin{equation}
\label{eq:gainfty2}
R_j^< \, \ge \, s_\gep \, \theta_j \sum_{n \le \gep/\theta_j} n K(n)
\exp(- \gd_j n) \stackrel{j \to \infty}\sim
s_\gep \gG(1-\ga) L(1/\gd_j) \, \left(\frac{\theta_j}{\gd_j}\right) \gd_j^\ga\,.
\end{equation}
To obtain \eqref{eq:gainfty2} we have used summation by parts, namely the
identity:
\begin{equation}
\sum_{n=1}^\infty nK(n) \exp(-\gd_j n)\, =\, 
\sum_{n=0}^\infty \Kbar (n) \exp(-\gd_j (n+1))\, -\,
\left( 1- \exp(-\gd_j)\right) \sum_{n=1}^\infty n \Kbar (n) \exp(-\gd_j n).
\end{equation}
On the other hand
\begin{equation}
\label{eq:gainfty3}
\left\vert R_j^> \right\vert 
\,\le\, \exp\left( -(\gd_j/\theta_j) \gep \right) \sum_{n> \gep/\theta_j}
K(n) 
\, \overset{j\to\infty}\sim \, \exp\left( -(\gd_j/\theta_j) \gep \right) \, \frac{L(1/\theta_j)}
{\ga}
(\theta_j / \gep)^\ga \,,
\end{equation}
therefore
\begin{equation}
\label{eq:gainfty4}
	\left\vert \frac{R^>_j}{R^<_j} \right\vert \,\le\,
	c \, \exp\left( -(\gd_j/\theta_j) \gep \right) \,
	\frac{L(1/\theta_j)}{L(1/\gd_j)} \, \left(\frac{\theta_j}{\gd_j}\right)^{\ga - 1}
	\,\le\, 	c' \, \exp\left( -(\gd_j/\theta_j) \gep \right) \,
	\left(\frac{\theta_j}{\gd_j}\right)^{\ga - 2}\,,
\end{equation}
where $c, c'$ are positive constants
 (we have explicitly used the fact that, for every $c_1>1$ and every $c_2>0$ there exists $c_3>0$ such that 
$L(x)/L(y) \le c_1 (x/y)^{c_2}$ whenever $x/y\ge c_3$
\cite[Th.~1.5.6]{cf:RegVar}). Therefore $|R^>_j / R^<_j| \to 0$ as
$j\to\infty$ and for $j$ sufficiently large we have
\begin{equation}
\sum_n K(n)
\exp(-\gd_j n) \sin (\theta_j n) \,\ge \, 
\frac 12 s_\gep \gG(1-\ga) L(1/\gd_j) \frac{\theta_j}{\gd_j} \gd_j^\ga\,,
\end{equation}
and then also in this regime \eqref{eq:aim} holds.

\medskip

The marginal case of $\ga =1$ and $\sum_n n K(n)=+\infty$ is treated as follows.

If $\ga\in [0, \infty)$ for the step analogous to \eqref{eq:Rieman2.1}
we split the sum according to whether $\theta_j n\le \gep$ or
$\theta_j n > \gep$. 
Summing by parts we obtain
\begin{equation}
\label{eq:byparts2}
\sum_{n=1}^N n K(n)\, =\, \sum_{n=0}^{N-1} \Kbar (n) \, -\, N \Kbar (N)
\stackrel{N \to \infty}\sim
\sum_{n=1}^N \frac{L(n)}{n} \, =:\, \hat L(N)
, 
\end{equation}
where in the asymptotic limit we have used \cite[Prop.~1.5.9a]{cf:RegVar}
that guarantees that $\hat L (\cdot)$ is slowly varying and 
that $\lim_{n \to \infty} \hat L(n)/L(n)=+\infty$.
From this we directly obtain that the first term in the splitting, {\sl i.e.} the sum over
$\theta_j n\le \gep$,
 is bounded below
by a positive constant, depending on $\gep$ and $\gamma$
(this constant can be chosen bounded away from zero for
$\gamma $ in any compact subset of $[0, \infty)$)
times $\theta_j \hat L( 1/\gd_j)$. The rest instead is bounded, in absolute value,
by a constant (independent of  $\gamma$)
times $\theta_j  L( 1/\theta_j)$, for $j$ sufficiently large  (just use
$\vert \sin (\theta_jn ) \exp(-\gd_j n) \vert \le 1$).
Using once again $\hat L(n) \gg L(n)$ for large $n$,
we obtain that $\sum_n K(n) \exp(-\gamma_j n) \sin (\theta_j n) >0$
for $j $ sufficiently large.

If instead $\gamma=+\infty$ we restart from \eqref{eq:splitnew}
and, by proceeding like in \eqref{eq:gainfty2} and  
\eqref{eq:gainfty3}, we obtain that for $j$
sufficiently large
\begin{equation}
\label{eq:splitnew2}
\sum_n K(n)
\exp(-\gd_j n) \sin (\theta_j n) \,\ge \, 
\frac 12 
s_\gep  \hat L(1/\gd_j) \, \left(\frac{\theta_j}{\gd_j}\right) \gd_j
\, -\, 2
 \exp\left( -(\gd_j/\theta_j) \gep \right) \, {L(1/\theta_j)}
\theta_j / \gep ,
\end{equation}
which is positive for $j$ sufficiently large
and the case $\ga=1$ and $\sum_n n K(n)=\infty$ is under control.

\medskip

Let us now consider the case of $\ga >1$, together with  
the case  $\ga=1$ and $\sum_n n K(n)<\infty$ and note that
in the latter case $L(\cdot)$ vanishes at infinity.  

In these cases for every $\gamma \in [0, \infty]$ we use the splitting in \eqref{eq:splitnew}
and for $j$ sufficiently large we have
\begin{equation}
\label{eq:Rieman2.1.2}
\sum_n K(n) 
\exp(-\gd_j n) \sin (\theta_j n)
\, \ge \, 
\frac 12 s_\gep \theta_j \sum_n n K(n) \, -\, \frac 2\ga L(1/\theta_j) \theta_j^\ga \gep^{-\ga},
\end{equation}
and the right-hand side is positive (again, for $j$ sufficiently large).
This concludes the proof of Proposition~\ref{th:Rb}.
\qed

\section{Sharp estimates}
\label{sec:sharp}

Throughout this section $K(\cdot)$ satisfies \eqref{eq:strongHyp}, we assume $b>0$  and we 
set $\nabla u_b(n):= u_b(n)-u_b(n-1)$ for $n=0,1, \ldots$ ($u_b(-1):= 0$).
We also introduce the discrete probability density $\mu_b$
on $\N \cup \{ 0\}$ defined by
\begin{equation}
\mu_b (n) \, := \, {\overline{K_b}(n)}/{m_b}, 
\end{equation}
with $m_b := \sum_n n K_b(n)$ and $ \overline{K_b}(n)
:= \sum_{j>n} K_b(j)$. Let us observe  that
\begin{equation}
\label{eq:Kbbar}
m_b \mu_b (n) \, =\, 
{K_b(n)} \sum_{j=1}^\infty \frac{K(n+j)}{K(n)} \exp(-bj) \stackrel{n \to \infty}\sim
\frac{1}{\exp(b)-1}K_b(n),
\end{equation}
and that this directly implies  the properties
\begin{equation}
\label{eq:alphabeta}
\frac{\sum_{j=0}^n \mu_b(j) \mu_b (n-j)}{\mu_b(n)} \stackrel{n \to \infty}\sim
 2 \hat {\mu_b} (\exp(b))
\ \ \ \textrm{ and } \ \ \
\frac{\mu_b(n+1)}{\mu_b(n)}\stackrel{n \to \infty} \sim \exp(-b). 
\end{equation}
We point out also that from \eqref{eq:rec} we get
\begin{equation}
\label{eq:grubel}
 \hat{\nabla u_n}(z) \, =\, \phi_b \left(\hat {\mu_b}(z)\right),\ \ \ \text{ with }
 \ \ \  \phi_b(z)\, := \, \frac 1{m_b z},
\end{equation}
at least for $\vert z \vert <1$,
like for \eqref{eq:BL}. 
 Of course the domain of analyticity of $\phi_b(\cdot)$ is 
 $\bbC\setminus \{0\}$ and if we observe that, by direct computation, we have
 \begin{equation}
 \hat {\mu_b}(z) \, =\,  \frac{1- \hat {K _b} (z)}{m_b(1-z)},
 \end{equation}
one can then extend the validity of \eqref{eq:grubel}
to all values of $z$ satisfying $\vert z\vert \le \exp(b)$ 
and $\vert z\vert<
  \inf\{\vert \zeta\vert>1 :\, \hat {K_b} (\zeta)=1\}$.

\medskip

\noindent
{\it Proof of Theorem~\ref{th:main}(2).}
Let us choose $b < b_0$.
We observe that the two properties in
\eqref{eq:alphabeta} are the hypotheses ($\ga$)
and ($\gb$) of \cite[Theorem~1]{cf:CNW}.
Hypothesis ($\gamma$) of the same theorem,
that is that $\hat{\mu_b}(z)$ converges at its radius of convergence
($\exp(b)$), is verified too.
Since $b<b_0$,  
 $\{\hat{\mu_b}(z): \, \vert z \vert \le \exp(b)\}\subset \bbC\setminus \{0\}$, 
  {\it i.e.} the range of the power series $\hat {\mu_b}(\cdot)$ is a subset of the  analyticity domain of $\phi_b (\cdot)$. 
Therefore  \cite[Theorem~1]{cf:CNW} yields
\begin{equation}
\label{eq:fromCNW}
\nabla u_b(n) \stackrel{n \to \infty}\sim \phi_b ^\prime \left( 
\hat{\mu_b} (\exp(b)) \right) \, \mu_b (n) \, =\, -\frac{\mu_b (n)}{\left(
 \hat{\mu_b} (\exp(b))
\right)^2 m_b },
\end{equation}  
and by \eqref{eq:Kbbar} we have
\begin{equation}
\nabla u_b(n) \stackrel{n \to \infty}\sim
-\frac{c(b) (\exp(b)-1)}{(c(b)-1)^2} K(n) \exp(-bn).
\end{equation}
We conclude by observing that
this yields
\begin{equation}
u_b(n) \, =\, - \sum_{j>n} \nabla u_b (j)  \stackrel{n \to \infty}\sim
\frac{c(b) }{(c(b)-1)^2} K(n) \exp(-bn)\, =\, \frac{K_b(n)}{(c(b)-1)^2},
\end{equation}
and the proof is complete. 
\qed

 \section{Some examples and further considerations}
 \label{sec:examples}

 Recall that $\gG (z):=\int_0^\infty t^{z-1} \exp(-t) \dd t$
for $\Re(z)>0$, that $\gG(\cdot)$ can be extended as a meromorphic function
to $\bbC$ and  that $\gG (z+1)=z\gG (z)$ for $z \notin \{0,-1,-2, \ldots\}$ (therefore
$\gG(n)=(n-1)!$ for $n \in \N$). 
 Much of the content of this section is based on the fact that
for $\gb \in \R \setminus\{0, -1,-2, \ldots\}$ and $\vert x\vert <1$ we have 
\begin{equation}
\label{eq:Gammaexp}
\sum_{n=0}^\infty \frac{\gG(\gb+n)}{n!} x^n \, =\, \gG(\gb) (1-x)^{-\gb}.
\end{equation}
This
is
just a matter of realizing that for $n \ge 1$ 
\begin{equation}
\frac{\dd ^n}
{\dd x^n} (1-x)^{-\gb} \, =\, \gb (\gb+1)\ldots (\gb+n-1) (1-x)^{-\gb-n},
\end{equation}
and the formula is the Taylor expansion in $x=0$.

Since $\sign(\gG(\gb))=(-1)^{\lceil \vert\gb\vert \rceil}$ for $\gb<0$
($\vert \gb \vert \notin \N$)
the first terms of the series in \eqref{eq:Gammaexp} have alternating 
signs, but for $n$ sufficiently large the sign stabilizes and, by Stirling's formula
\begin{equation}
\label{eq:Stirling}
\gG(x)\stackrel{x \to \infty}\sim \exp(-x)x^{x-(1/2)}\sqrt{2\pi},
\end{equation} 
one readily sees that
$\gG(n-\ga)/n! \stackrel{n\to \infty} \sim  1/n^{1+\ga}$.
Therefore, with the help of 
\eqref{eq:Gammaexp} we
can build probability inter-arrival distributions with the type of decay we are interested 
in and for which the $z$-transform is explicit. 

\medskip

\begin{rem}\rm
It is not difficult to see that one can differentiate, say $j$ times, 
the expression in \eqref{eq:Gammaexp}
generating thus sequences which decay like $(\log n)^j/n^{1+\ga}$
and that, for sufficiently large $n$, do not change sign.
This provides examples involving slowly varying functions.
\end{rem}
\medskip

Since we are just developing examples and that generalizations
are straightforward, we
 specialize to the case of $-\gb=\ga\in (0,1)$.
 
\subsection{The basic example} 
\label{sec:basic}  

In this section we study the case of
\begin{equation}
\label{eq:basic}
K(n) \, :=\, 
\frac{\gG (n-\ga)}{-\gG(-\ga) \, n!} \stackrel{n \to \infty} \sim
\frac
{n^{-1-\ga}} {-\gG(-\ga)}.
\end{equation}
Note that $\sum_{n=1}^\infty K(n)=1$ follows from \eqref{eq:Gammaexp},
 with $\gb =-\ga$, as well as, with reference to  \eqref{eq:Kb},
$c(b)= 1/(1-(1-\exp(-b))^\ga)$ and
\begin{equation}
\label{eq:Khat}
\hat{K_b} (z) \, =\,\frac{ \left(1-(1-z\exp(-b)) ^\ga \right)}
{\left(1-(1-\exp(-b)) ^\ga \right)}.
\end{equation}
In defining  $z^\ga$ for $\ga$ non integer, 
we  choose the  cut line 
$\{z\in \R: \, z< 0 \}$.
With this choice $ (1-z\exp(-b)) ^\ga$, and therefore 
$\hat{K_b}(\cdot)$, has a discontinuity
 on  the line $\{z\in \R:\, z>\exp(b)\}$.
 
We observe  that, for every $b>0$, $\hat{K_b}(z)=1$ for $\vert z \vert \le \exp(b)$
only if $z=1$, therefore  Theorem~\ref{th:main} holds with $b_0=\infty$.

\medskip

\begin{rem}
\rm
In the special case under consideration, but also
in all the other cases considered in this section,
one can obtain and go beyond   Theorem~\ref{th:main} by direct computations.
In fact if we set $q(z):= (1-z\exp(-b))^\ga$ we have for $\vert q(z)\vert< \vert q(1)\vert$
\begin{equation}
\frac1{1-\hat {K_b} (z)}\, =\, \frac{1-q(1)}{q(z)-q(1)}\, =\, -\frac{1-q(1)}{q(1)} \, \sum_{j=0}^\infty
\left(
\frac{q(z)}{q(1)}
\right)^j.
\end{equation}
Now we set
\begin{equation}
\label{eq:set4.2}
R_m(z)\, := \, \gD_b(z)
\, +\,
\frac{1-q(1)}{q(1)} \, \sum_{j=1}^m
\left(
\frac{q(z)}{q(1)}
\right)^j,
\end{equation}
and we note that
$(q(z))^j=  (1-z\exp(-b))^{j\ga}$ and therefore once again
\eqref{eq:Gammaexp} provides the expansion for $(q(z))^j$
if $j\ga \notin \N$ and the $n$-th term in the 
power series (of $(q(z))^j$) behaves, as $n \to \infty$, like
$c \exp(-nb)n^{-1-j\ga}$, $c \neq 0$.
Note that if $j\ga \in \N$ the arising expression is just a polynomial
and hence does not contribute to the asymptotic behavior of the
series expansion. 

Finally, the series expansion $\sum_n r^{(m)}(n) z^n$ of $R_m(\cdot)$
 can be controlled by observing that this function is analytic in
 the centered ball of radius $\exp(b)$ and by using 
the formula
\begin{equation}
\label{eq:path}
r^{(m)}(n)\, =\, \frac 1{2\pi i}
\oint \frac{ R_m(z)}{z^{n+1}} \dd z\, =\, 
\frac {\exp(-bn)}{2\pi } \int_0^{2\pi}
{R_m\left( \exp(b+i\theta)\right)}
\exp\left(-i n \theta\right)
\dd \theta,
\end{equation}
where the contour in the middle term is (say) $\vert z\vert =r$, for $r\in (0,\exp(\gb))$,
and the last term is obtained by letting $ r\nearrow \exp(b)$,
using the fact that $R_m(\exp(b+i\theta))$ is bounded.
 In fact, from the explicit expression and by construction, one readily sees that
$R_m\left( \exp(b+i\theta)\right)$ is smooth except at $\theta=2\pi k$, $k\in \Z$,
where it is $C^{\lfloor (m+1) \ga \rfloor}$. By using the fact
that $n$-th Fourier coefficient of 
a $C^k$ function is $o(n^{-k})$, we see that
$r^{(m)}(n) \,=\, \exp(-bn)  o( 1/n^{\lfloor (m+1) \ga \rfloor})$.

The chain of considerations we have just made leads to an explicit  expansion to all orders for 
$\exp(bn)(u_b(n)-u_b(\infty))$ 
as a sum of terms of the form
$c_{j_1, j_2} n^{-j_1-\ga j_2}$, for suitable (explicit) real coefficients
$c_{j_1, j_2}$ ($j_1$, $j_2\in \N$). 
\end{rem}

\subsection{Singularities and slower decay of correlations}

From the basic example one can actually build a large number of
{\sl exactly solvable} cases that display the more general phenomenology 
hinted by Theorem~\ref{th:main}: in particular that, in general, 
 $b_0< \infty$.

For example, fix $m\in \N$ and define
\begin{equation}
\label{eq:basm}
K(n)\, :=\, \begin{cases}
{\gG (n-m-\ga)}/\left({-\gG(-\ga) \, (n-m)!}\right)
&\text{ for } n =m+1, m+2, \ldots
\\
0 &\text{ for } n =1, 2, \ldots, m.
\end{cases}
\end{equation}
Note that this is nothing but the previous choice of $K(\cdot)$
{\sl translated} to the right of $m$ steps. 
 Therefore
 \begin{equation}
\label{eq:Khat2}
\hat {K_b} (z) \, =\,z^m \, \frac{ \left(1-(1-z\exp(-b)) ^\ga \right)}
{\left(1-(1-\exp(-b)) ^\ga \right)}.
\end{equation}
 Once again the radius of convergence is $\exp(b)$, but this
 time, in general, it is no longer true that one cannot find  a solution $z_0$
 to $\hat{K_b}(z_0)=1$ in the annulus $1<\vert z_0\vert < \exp(b)$.
 
 \smallskip
 
 Let us choose $\ga =1/2$ and let us first look at the case of $m=1$.
 One can directly verify that 
 \begin{equation}
 z_0\, =\,
 -\frac 12 \left(1+
 \sqrt{8 \exp(b)  \left(1-\sqrt{1-\exp (-b)} \right) -3}
 \right) \, <\, -1, 
 \end{equation}
solves $\hat{K_b}(z_0)=1$, that it is the unique solution (except the trivial solution $z_0=1$),
 and $\vert z_0\vert < \exp(b)$ for $b>b_0$ with
\begin{equation}
b_0 \, := \,
\log \left(3/2+\sqrt{2}- \sqrt{\sqrt{2}+5/4 }\right)\, =\, 0.248399...
\end{equation}
So, if $b>b_0$, since $ z_0$ is a (simple) pole singularity
of $\gD_b(\cdot)$ we can write 
\begin{equation}
\label{eq:pole_sing1}
\gD_b(z) \, =\, \frac 1{z_0 K^\prime_b (z_0)\left(1-(z/z_0)\right)}\, +\, f(z),
\end{equation}
with $f(\cdot)$ a function which is analytic on the centered ball
of radius $\exp(b)$. Therefore
\begin{equation}
\label{eq:pole_sing2}
u_b(n) -u_b(\infty) \, =\, 
 \frac 1{z_0 K^\prime_b (z_0)} z_0^{-n} + \gep(n), 
\end{equation}
and $\limsup_{n \to \infty}(1/n)\log \vert \gep (n)\vert = -b$.
\smallskip
  
\begin{rem}\rm
Note that $z_0= -1-\exp(-b)/4+O(\exp(-2b))$ for $b $ large, so that 
the rate of converge of $u_b(n)-u_\infty (n)$  becomes
smaller and smaller as $b$ becomes large. 
\end{rem}

\medskip

Going back to \eqref{eq:basm}, for $m$ larger than $3$
one can no longer explicitly find all the solutions $z$ to $\hat{K _b}(z)=1$.
However we have the following:

\medskip

\begin{proposition}
\label{th:m}
For every $b>0$ and $\ga \in (0,1)$ one can find $m\in \N$ such that
if $K(\cdot)$ is given by \eqref{eq:basm}  then
there exists a solution $z_0$ to  $ \hat{K_b}(z_0)=1$ 
with $1<\vert z_0\vert< \exp(b)$.
\end{proposition}

\medskip

\begin{rem}\rm 
In general, once the solutions to $ \hat{K_b}(\cdot)=1$ of minimal absolute value
(in the annulus  $\{z:\, 1<\vert z\vert<\exp(b)$) are known, it 
is straightforward to write the sharp asymptotic behavior
of $u_b(n)-u_b(\infty)$. For example if $z_0$ is a complex solution,
then also its conjugate is a solution. If these have minimal absolute value
among the solutions and if they are simple solutions, 
for a suitable (and computable) real constants $c_1$ and $c_2$ ($\vert c_1\vert+
\vert c_2\vert >0$)
  we have
\begin{equation}
\label{eq:m}
u_b(n)-u_b(\infty) \stackrel{n \to \infty}\sim  \vert z_0\vert^{-n} 
 \left( c_1 
 \cos\left( n \arg\left(z_0\right)\right)+ c_2 \sin\left( n \arg\left(z_0\right)\right)
 \right).  
\end{equation}
 An analogous formula is easily written 
 in the general case.
\end{rem}

\medskip

\noindent
{\it Proof of Proposition~\ref{th:m}.}
In reality, we are going to do something rather cheap, but we are actually proving
more than what is stated: we are going to show that for every $b>0$
and every $r\in (0, \exp(b))$ we can find an $m$
such that there are $m$ zeros of $\hat{K_b} (\cdot) -1$
in the annulus $\{z:\, 1<\vert z\vert < r \}$.

Given  $b>0$, since the only solution $z$ to
$1-(1-z \exp(-b))^\ga)=0$ is $z=0$,
then for every $r\in (1, \exp(b))$ we have
\begin{equation}
x_r:=
\inf_\theta \left\vert\frac{1-(1-r \exp(-b+i\theta))^\ga)}{1-(1- \exp(-b))^\ga)}\right\vert\, >\, 0.
\end{equation}
Therefore (recall \eqref{eq:Khat2})
 $\vert \hat {K_b} (z)\vert \ge r^m x_r$, if $\vert z \vert =r$.
 Therefore for $m$ sufficiently large we have 
  $\vert \hat {K_b} (z)\vert >1$ for $\vert z \vert =r$: let us fix such a couple  $(m,r)$.
Rouch\'e's Theorem ({\sl e.g.} \cite[p.~153]{cf:Ahlfors})  guarantees that if $f$ and $g$ are analytic in a simply connected domain
containing the simple closed curve $\gamma$ and if 
 $\vert f(z) -g(z)\vert < \vert f(z)\vert$ for $z \in \gamma$, then
 $f$ and $g$ have the same number of zeros enclosed by $\gamma$.
 Let us apply Rouch\'e's Theorem  with
 $f(z):= \hat{K_b} (z)$ and $g(z):= 1-\hat{K_b} (z)$ and 
  $\gamma:= \{z: \, \vert z \vert =r \}$, so that
 $\vert f(z) -g(z)\vert =1 <  \vert f(z)\vert$ for $z\in \gamma$, by the choice of $m$.
 But $ \hat{K_b}(\cdot)$ has precisely $m+1$ zeros (they are all in  $0$)
 and therefore also $1-\hat{K_b}(\cdot)$ has $m+1$ zeros enclosed by $\gamma$. 
 Of course $1-\hat{K_b}(\cdot)$ has a zero in $1$ and all the other
 zeros have absolute value in $(1, r)$. 
 \qed

\section*{Acknowledgments}
 I am greatly indebted with Bernard Derrida for having supplied
the basic example of Section~\ref{sec:examples} and for several
discussions. I am also very grateful to Francesco Caravenna and
to Fabio Toninelli for important observations and discussions. 
The author acknowledges the support of ANR, project POLINTBIO.

\end{document}